# HISTORY OF MATHEMATICS EDUCATION - ITALY

**Marta Menghini**
*Dipartimento di Matematica, Università di Roma Sapienza, Roma, Italy.*

**Keywords:** Mathematics education, history of education, Italy, school-mathematics, teaching of geometry, teaching of arithmetic.

**Contents**

1. Introduction: System of education in the period of the unification of Italy
2. The role of mathematics in the Italian educational system
3. Mathematics in primary school and the training of the teachers
4. Main features of Italian mathematical instruction
5. The reform of 1923 and the birth of the scientific Lycée
6. Attempts for a reform of secondary schools: 1950 to 2010
   Bibliography

**Summary**

Starting from the Italian unification in 1861, we describe the role of mathematics teaching in the Italian system of education giving an account of the main features of Italian mathematics education, and referring to the official syllabi, to the used textbooks, to the debates. We will see that interesting proposals in a period of about 150 years will always have to face a tendency to rigor as well as a resistance towards innovation.

**1. Introduction: System of education in the period of the unification of Italy.**

The first law to regulate Italian schools was the Casati Law of 1859 (R.d. 13.11.1859, n. 3725). It was originally passed only for the Kingdom of Sardinia (including Piedmont) and for Lombardy, and was later gradually extended to the other Italian regions after their annexation in 1861 and 1871. The Casati Law established the general characteristics of state education.

As to secondary education, the law distinguished between a *classical education*, whose purpose was the acquisition of that literary and philosophical culture which opened the way to the State Universities (Art. 188), and a *technical education,* which aimed at providing a general culture to young people intending a job in the public services, in industry, in commerce and in agricultural management (Art. 272). Secondary education consisted of a first and a second level: to cover classical secondary education, the Casati Law introduced the Gymnasium and the Lycée (*Ginnasio-Liceo*), which were to become the point of reference for all Italian secondary education (Vita, 1986). The Technical School (*Scuola Tecnica*) and the Technical Institute (*Istituto Tecnico*) were set up for technical secondary education (Table 1.) Pupils entered the Gymnasium and the Technical School after a primary school that extended over four years. The Technical School thus covered the age range of the present-day middle school (11-14), while the Gymnasium lasted for five years and hence included the first two years of high school. The Technical School soon lost its characteristic of being a preparatory school for the Technical Institute and was transformed into a school for general education. After the Gymnasium, pupils completed their classical education by attending the Lycée for three years. Initially, a technical education was shorter than its classical counterpart because, after the three years of Technical School, only a further three years were foreseen at the Technical Institute. However, through an 1871 reform (circ. 17.10.1871) the duration of the Technical Institute was extended to four years and, in some cases, to five.

In 1896 a Complementary Course (parallel to Technical School and to the *lower Gymnasium*, i.e. the first three years of the Gymnasium) was introduced – followed by the *Scuola Normale* - specifically for the *Istruzione Magistrale*, which was the cycle aimed at educating primary school teachers. The norms related to *Istruzione Magistrale* were included among the norms for primary instruction.

Table 1. Instruction in Italy at the end of 19$^{th}$ century

**2. The role of mathematics in the Italian educational system**

The period of the unification was without a doubt a unique period in the scientific history of Italy.

Patriotism intertwined with a lay and positivist mentality born with the *Risorgimento* (the Italian Resurgence), which was widespread amongst the bourgeoisie, gave scientific research a central role that it was not to have again (Brigaglia & Masotto, 1982). Mathematical research flourished and mathematics teaching was considered relevant at all levels of education.

**2.1. Mathematics in classical and technical instruction**

The characteristics of the teaching of mathematics in classical and technical instruction are well defined by the Law of Education Minister Coppino in 1867 (R.d. 10.10.1867). As for the classical instruction, we read, in the first page, that

*Mathematics in classical secondary schools is not to be regarded only as a set of propositions or theories, which are useful in themselves [...] to be applied to the needs of life; but primarily as a means of intellectual culture, as a gymnastics of the mind aimed to develop the power of reasoning, and to help that right and sane criterion needed to distinguish what is true...*[all translations are by the author].

The teaching of mathematics in the Gymnasium-Lycée covered initially three classes only (the second year of the upper gymnasium and the first two years of the Lycée) with a high number of hours per week (5 + 6 + 7½), then extended to all the five years of high school (4 + 4 + 7½ + 4½ + 2). The programmes of 1867 reduced the previous mathematical contents in order to focus on the rational development of geometry and arithmetic (see 4.1 and 4.2), but the high number of hours devoted to mathematics testifies for the credit deserved in that period to this subject. In the next decades the number of hours of mathematics in the Lycée will harshly decrease (to 2+2+3+3+2 hours per week) and the abolishment, at the end of the 19$^{th}$ century, of the written examination in mathematics at the end of the Gymnasium and of the Lycée shows the diminished value ascribed to mathematics.

According to the same Coppino Law of 1867, mathematics had a different role in technical education, in that:

*The aim of the teaching of mathematics in technical schools is to provide young men, in a limited time, with the greatest amount of knowledge useful for applications in the arts and crafts.*

The Technical Institutes had developed considerably after the unification of Italy, having undergone various reforms since 1860, all of which recognised the necessity of separate developments of humanist and technical education, with an eye to the model of the German Realschulen (Morpurgo, 1875, XXVI and on; Ulivi, 1978). At the beginning of the 20$^{th}$ century students in the technical instruction were about 60% of all secondary school students (Scorza, 1911). A reform of 1871 recognised the need for a general literary and scientific education in technical education and instituted a physics-mathematics section (*sezione fisico – matematica*). Since this section did not have the aim of qualifying students to go into the professions, and permitted university entrance (to the Faculty for Engineering), it could be seen as the scientific alternative to the Lycée. The first two years of the Technical Institute were common to all the sections and had a general and preparatory character. In the second "advanced" two-year phase, the syllabus of the third year was common to the physics-mathematics section and to the industrial section. In the fourth year the program was specific for the physics - mathematics section. Mathematics had 6+5 hours per week in the first two years, and 5 + 5 in the two further years of the physics-mathematics section.

**2.2. Arenas for a debate**

After the unification, Italian mathematics teachers acquired a professionalization evidenced by the foundation of journals addressed to them. The first one, *Rivista di Matematica Elementare*, appeared in 1874 (Furinghetti & Somaglia, 1992). Before this year, only a few articles were published in journals addressed at professional mathematicians and students of University. This was the case of *Giornale di Matematiche,* which hosted the debate about the introduction of Euclid's Elements in the Gymnasium (see 4.1.). An important journal in the field of mathematics education was *il Periodico di Matematica,* which appeared in 1886. At the beginning of the 20th century, the *Periodico* will share the stalls of mathematics teachers and educators with another journal, the *Bollettino di Matematica,* which in 1949 will change its name in *Archimede*. Both journals still survive.

In 1895 a group of mathematics teachers founded the association Mathesis (Giacardi, 2004). The aim of this association was the improvement of the scientific and didactical knowledge of the teachers. These aims had to be pursued by facing the problems concerned with school syllabi, methods of teaching, teacher training, schoolbooks, etc. The first president was the teacher Rodolfo Bettazzi. This was a period of a big interchange between school and university. May be due to the studies on the foundations of mathematics, many mathematicians showed an interest in school mathematics, not only by writing or translating manuals for the school but also by intervening in the political and cultural discussions, in teacher training, in the construction of curricula. On the other hand good teachers, as Bettazzi himself, hold lectures at the university. The association Mathesis grew rapidly and included among its members not only teachers but also university professors. Its members participated in debates about the teaching of mathematics in school, expressing their opinions in the Congresses of the Association and on its Bulletin (*Bollettino dell'Associazione Mathesis*). We will meet them again in the next sections.

Among the presidents of Mathesis we find, in 1910, Guido Castelnuovo. As a member of the Italian subcommittee of the International Commission on Mathematical Instruction (ICMI), Castelnuovo established important international contacts. The *Bollettino della "Mathesis" Società Italiana di Matematica,* published as an autonomous and rich journal in 1909-1920, hosted translations of international debates and the reports of the Italian subcommittee. In this way Mathesis participated in the discussions about the introduction of calculus in schools and about the applications of mathematics; many mathematicians, even within Mathesis, opposed the introduction of calculus as well as the utilitarian aspect of mathematics (see Fehr, 1911). When the *Bollettino* ceased the publication, *il Periodico di Matematica* became the official organ of Mathesis. Federigo Enriques was another well-known mathematician who became president of Mathesis (in 1919). Enriques succeeded in enlarging the Association such as to have a major incidence on the governments. Mathesis will in fact become a counterpart in important questions, like the problem of the standardization of the curricula in the regions that formerly belonged to the Austrian Empire, as Trento and Trieste (Zuccheri & Zudini, 2007).

**2.3. Towards a modern Lycée**

Even if the physics-mathematics section of the Technical Institute could represent an alternative to the classical curriculum, the international movement towards new scientific humanities was felt also in the Italian milieu. The first signal in this direction can be noted already in 1875, when a decree of the Minister Bonghi (R.d. 7.1.1875, n. 2337) allowed students ending the Lycée to pass the school-leaving examination even if they, having a good note in mathematics, failed in Greek (or failed in mathematics having a good note in Greek). In 1899 a decree of Minister Baccelli (d.m. 3.11.1899) allowed the students of some particular Lycée of bigger towns to choose, in the last year, between Greek and mathematics. New topics were added for those who chose mathematics, as the harmonic division, homothety, spherical triangles and the calculation of approximate values of $\pi$. This possibility was suppressed the next year, but in 1904 a new reform by Minister Orlando (R.d. 11.11.1904) again allowed to choose, in the last two years of the Classical Lycée, between Greek and mathematics. This reform was a first step towards a modern lycée, but mathematicians themselves opposed it as contrasting the educational value of the Lycée (*Bollettino della Mathesis*, 1911).

In 1905 the Minister of Education Bianchi appointed a Royal Commission for secondary education to explore the possibility of creating a unified middle school, followed by an upper education divided into three branches. The Commission had a tormented life. Different views, including the difficulty of accepting the idea of a complete equalization of the classical and the scientific studies, led to the resignation of a part of its members and to the publication of their opinions in a volume (Galletti & Salvemini, 1908). The Commission, which had continued the work, published its final report in 1909 (Commissione Reale, 1909). Giovanni Vailati edited the part concerning mathematics. The Commission foresaw, inter alia, three Lycées: a classic, a modern and a scientific one. The programmes of the *Commissione Reale* never entered into force. It is worth to mention some points concerning the scientific Lycée. In the first year, geometry is dealt with only through problems of construction, while in the second year it is introduced in a more systematic way. The third year requires the "graphical representation of the variance of the expressions of the second degree", not re-appeared until the 1980s.The concept of derivative is introduced, and its different aspects are emphasized such as the "rate of growth of a function in a given interval or in a given point." In the fifth year we find probability, besides integral calculus. Note that the "representation of statistical data" appears only in the modern Lycée. Assuming that the learning process goes from concreteness to abstraction, the comments state that pupils should not learn "theories before knowing the facts to which they refer, and from which their meaning can be obtained by abstraction".

*2.3.1. The Modern Lycée: 1911 – 1923*

The *Liceo Moderno* was the only partial implementation of the project presented by the *Commissione Reale*. It was founded in 1911, only in some provinces, as a section of the Lycée. The aim to prepare students for (not necessarily scientific) university studies is achieved through the study of Latin, modern languages and sciences. Mathematics is important as a language suitable to describe natural phenomena and therefore its contents, edited by Castelnuovo (1913), show significant novelties. "The renewal of mathematics in the seventeenth century is linked to the flourishing of experimental sciences. In this view the teacher will point out that the fundamental concepts of modern mathematics, that of function in particular, are suggested by observational sciences, which then exert a beneficial influence on the development of mathematics itself". The notion of a function and infinitesimal calculus were introduced, in an official way, for the first time in this school, with the suggestion to use an experimental and inductive approach - which completes the deductive method - and to harmonize it with the course of physics.

*2.3.2. The programmes of the Ispettorato Centrale dell'Istruzione Media*
A few years later, the Ministry of Education (1918) proposed a reorganization of the curricula of all secondary schools, edited by the general director Vittorio Fiorini, former member of the *Commissione Reale*. No changes were suggested for the *Liceo Moderno*, but the proposals of the *Commissione Reale* for the physics-mathematics section of the Technical Institute were followed. The programmes took note of the fact that in many of these schools basic notions of calculus were already taught. It is recommended to introduce these topics by describing the historical process and the link with the sciences that created them. The use of derivatives in maximum problems and of the integrals in the calculation of areas and volumes will allow solving problems that previously required complicated algebraic transformations. These programmes, never adopted officially, represent the most advanced vision of the role of mathematics in the school when the Gentile Reform entered into force in 1923.

## 3. Mathematics in primary school and the training of the teachers

## 3.1. The Elementary School

The Casati Law also set up Elementary School. It lasted 4 years, only the two first years (*inferior course*) being compulsory. Local institutions paid for the teachers. In 1904 Elementary School was extended to two three-year courses. Pupils who wanted to enter secondary instruction could leave Elementary School after the fourth year, with a final examination. Compulsory Elementary School was extended, but depending on the number of classes existing in the municipality. Pupils could not enter Elementary School before their sixth year, and in any case had to leave the second course when becoming fifteen.

The teaching of mathematics in Elementary School was mainly devoted to arithmetic. One couldn't find geometry in the first three years, even if geometrical contents could be found in the teaching of fractions and of weights and measures (see Conti, 1911). In the second three-year course geometry teaching dealt with shapes, names and rules for the measure of the simplest geometrical figures.

It is interesting to note that teaching in Elementary Schools was influenced by psychologists and by educational sciences; this didn't happen in secondary school. We find, for instance, the interesting proposals of Anton Maria Bustelli (1889), who suggested – in the second course (two-year course at that time) - three phases for each new geometric form: 1. *presentation* of the object or of the figure (a drawing, cutted paper, pieces of wood, ...); 2. *description*, which means to undo and to redo the figure analyzing the relations between its parts; 3. *definition* of the figure. Maria Montessori will take up these phases again in her geometric work (Montessori 1934); and similar phases will be described in the work of Pierre van Hiele (1958).

The programmes established in 1905 (R. d. 29.1.1905) had a practical character. They suggest to avoid an excessive amount of written calculations (with respect, e.g., to the programmes of 1888) and to apply mathematics to the needs of *domestic* life or of little trades. An intuitive – experimental method is suggested as well as the cyclic development of the topics; this means that the four operations are all treated in the first class with numbers from 1 to 20, then in the second class with numbers from 1 to 100, etc. Moreover we find the suggestion to avoid technical terms (factor, quotient) and not to overcharge the children's minds. The teaching of fractions must not go beyond the needs of a popular elementary school. In the last classes calculations must be proposed in the form of problems. Examples from real life have to be "effectively real". This is particularly important in the last year, when the programmes include monetary calculations (as interest, discount). The school had to furnish the class with didactic materials as: a collection of common weights and measures; a cube that can be decomposed into 8 little cubes, a cylinder, a cone, a sphere and a pyramid made of cardboard or wood; of course other materials could be added, in particular those for the teaching of geometry and for an experimental justification of rules for measure.

Starting from the Gentile reform of 1923 (see 5.1.), Elementary School will last five years. This reform stressed the literary and artistic aspect in education. The teacher had a certain freedom, but the programmes were written under the influence of Fascist culture. The teaching of mathematics included only arithmetic.

In 1945 a Commission was set up by Italian Minister De Ruggiero with the collaboration of the Allied Countries, and chaired by the pedagogist Colonel Carlton Washburne, a follower of John Dewey. New curricula for Elementary School and Kindergarten were launched (d.m. l9.2.45 n.459 and d.l.gt 24.5.45 n.549). The teaching of arithmetic and geometry had to take into account the mental images and the intuitions of the pupil. The contents and the instructions looked back at the "pedagogical" programmes of the beginning of the century; the difficult political and economic situation in Italy hindered their success.

**3.2. The Normal School**

Teachers of Elementary School were formed in the Normal Schools. These lasted 3 years, and were divided into Normal Schools for males (to teach in the inferior elementary courses for males) and Normal Schools for females (to teach in female or mixed courses, and in all superior courses). Only at the end of the nineteenth century a specific Complementary three year course (*Corso Complementare*), was set up to precede the Normal School.

The mathematical contents of the Normal School were approximately those of the Elementary School, at a not too higher level. And, of course, those of the Corso Complementare were at an intermediate level, including – from time to time – rational arithmetic. In 1867 we find in the Normal School practical arithmetic with elements of book-keeping (Di Sieno, 2008). Geometry consists in the calculation of areas and volumes; in applications to measurement of land using the most common geodetic instruments (particularly in the agricultural communities); in the conversion from old measures to new ones and vice versa (due to the unification of the metric system). Before 1883 we find in the programmes integer and decimal numbers, the decimal metric system, fractions. In 1890 rational arithmetic was added to practical arithmetic. In 1892 the contents again include only practical arithmetic, while plane geometry seems to be some more extended. The becoming teacher had to be trained in the "exactness of the language and of the signs", and had to take his examples from domestic life, physics, geography and the other sciences.

**3.3. The Kindergarten and the influence of educational sciences**

To each Normal School a Kindergarten (*Giardino d'infanzia*) was annexed. There weren't programmes for that kind of school, but the curricular planning often mentions mathematical topics.

At the beginning of the 20$^{th}$ century teaching was dominated by the ideas of Friedrich Fröbel, who used concrete materials (gifts), to bring the pupils to discover their properties by playing. The "gifts" were among the materials that a school had to supply. Follower of Dewey's idea of active learning, M. Montessori created a school for children called *Casa dei bambini* (Children's home) in Rome in 1907; this school, apart from being a place of didactical research - often devoted to children with particular needs - influenced notably the future teaching in the lower grades. The main mathematical content of these experimental teaching was arithmetic, but the simple fact of working with real objects, of observing space and forms, of drawing, of handling concrete materials allowed an approach to geometry. Montessori assumed that manipulatives allow the discovering of relations between the various parts of the figures. Besides the well-known pink tower (cubes with different sides to be ordered to form a tower), Montessori used geometric forms that had to fit in the holes of a wooden tablet and iron geometric forms that could be used as help for drawing. Though being already active at the beginning of the 20$^{th}$ century, Montessori influenced teaching in Kindergarten and Elementary School only some decades later.

**3.4. Developments in the training of teachers**

For a long time the teachers of Elementary Schools were the only ones to receive a specific training in the Normal Schools. In 1923 these became *Istituti Magistrali*, and the duration became of four years.

To teach in secondary school a University degree (usually in mathematics) was needed, without particular courses devoted to didactical or pedagogical issues. The question of a specific training for pre-service teachers had often been discussed, but without effects. A Royal decree of 17.8.1862 confirmed the *Scuola Normale* of Pisa as the only Italian institution apt to educate secondary teachers. As it remained in fact the only one, addressed to few selected people, it slowly became a high level scientific institution, devoted to the training of University professors rather than of school professors. In 1875, a decree of Minister Bonghi (22.10.1875, n. 247, Art. 23) allowed some University Faculties of Science to give a special diploma for teachers, to be obtained with courses of a school (*Scuola di Magistero*), additional to the courses of the degree (in Mathematics, in our case). These schools too didn't bring to an establishment of teacher training and had frequent rearrangements (see http://math.unipa.it/~grim/storia_didmat_it/) till to their suppression in 1920.

Among the few successful experiences there was that of Turin at the beginning of the 20$^{th}$ century, where well-known mathematicians interested in teaching, as Corrado Segre and Enrico D'Ovidio, collaborated (the evolution of the teaching profession in Italy, in pre- and post-unity Italy, is illustrated by Furinghetti & Giacardi, 2012).

In 1921 Minister Corbino set up a degree in Mathematics and Physics, devoted to those who wanted to become teachers (but not compulsory to become a teacher). This degree included a new course of Complementary Mathematics (*Matematiche complementari*) and a course of Physics Laboratory, but no pedagogical subjects. In 1923, the Gentile Reform established that mathematics and physics in high school had to be thought by the same teacher. The mixed degree in mathematics and physics was suppressed in 1960; the course of *Matematiche Complementari* remained as a stable and typical course in the curriculum for the degree in mathematics with an orientation towards teaching. The suppression of the mixed degree didn't reach the aim of separating the teaching of mathematics and physics, as many teachers hoped. On the contrary, the middle school reform of 1962 presented the unified teaching of "mathematics and scientific observations".

Only in 1990 a law (19.11.1990, n.341) will create a degree for Elementary School teachers (*Corso di Laurea in Scienze della Formazione Primaria)*, and establish a post-degree two-year specialized school for secondary teachers. The SSIS (*Scuola di specializzazione all'Insegnamento Secondario*) will start only in 1999. It will include pedagogical subjects, history (of mathematics), laboratories, teaching practice and work experience in school. The SSIS was suspended (law n.133, 2008) and substituted by a one-year TFA (*Tirocinio Formativo Attivo,* d.m. 10.09.2010), not yet activated (in 2012).

## 4. Main features in Italian mathematical education

### 4.1 Rational geometry in the Lycée from 1867

*4.1.1. The debate*

The law by Minister Coppino in 1867 (see 2.1.) introduced the use of Euclid's Elements as the textbook for the teaching of geometry in the Gymnasium-Lycée. More specifically, Book I was intended for the fifth year of the Gymnasium (15 years old pupils), Books II and III for the first year of the Lycée, and successive Books for the second year of the Lycée. The reform reduced the previous contents of mathematics, which included combinatorics, polarity, conic sections, inequalities, maximum and minimum problems, limits. We will find these topics later on only in the physics-mathematics section of the Technical Institute.

The mathematician Luigi Cremona was a member of the commission that brought Euclid's text into the schools. He also contributed in an unofficial way to the Italian edition of Euclid's Elements that was edited by Enrico Betti and Francesco Brioschi (Betti & Brioschi, 1867), both consultants for the Ministry of Education. The reform sparked off some criticism, but it has provided the model for the teaching of geometry right up to the present day (Maraschini & Menghini, 1992; Menghini, 1996). Geometry was seen as *mental gymnastics* and was intended to get the young used to the rigour of reasoning: only the pure geometry of the Elements could carry out such a task. In truth it is not difficult to interpret the 1867 reform as a desire to create a new Italian model freed from the foreign textbooks that were then common in Italy. This is also highlighted in a letter of Cremona to Betti of 8.9.1869 (Gatto, 1996):

*People can say what they want but Euclid's is still the most logical and most rigorous system we have: all the successive systems are impure hybrids; in seeking to remove one defect, they fall into other worse ones and more than anything else they stop being true geometric systems. Legendre suffices as an example and this even though he is the most respectable of the elementary geometry reformists. However, if we cast our minds back to the books used in our schools before 1867, which would now be re-introduced if the syllabuses were modified, who would dare to deny that the introduction of the Euclidean method has been of immense benefit for our schools?*

The structure of Adrien-Marie Legendre's book was very similar to the Elements of Euclid, but Cremona criticized the fact that some proofs were simplified and the author sometimes made use of arithmetic and algebraic notations that hid the pureness of the geometric treatment. Some years later Felix Klein was to observe (Klein, 1909):

*The main input to mathematics teaching came from L. Cremona, known in the scientific field for his development of modern geometry. [...] for school purposes he advised a logical approach, whilst in his teaching and scientific work he has always stressed intuition. The link between these two visions of Cremona is not clear.*

We must anyway remember that for technical instruction Cremona suggested a different approach to geometry, based primarily on applications of projective geometry (see section 4.2). The disagreement of a part of the Italian mathematical community, generated by the introduction of the Elements, is well documented in several articles published in the *Giornale di matematiche* run by Giuseppe Battaglini (Hirst, 1868; Wilson, 1868; Brioschi & Cremona 1869). The articles by James M. Wilson and Thomas A. Hirst (who had a long correspondence with Cremona between 1864 and 1880) are translations of two articles which appeared in 1868 in the British magazine *The Educational Times*, and do not refer directly to the Italian situation, but rather to a similar and much older problem in England. They were inspired by new geometry textbooks for British schoolchildren that had just come out and represented two attacks against Euclid. There are criticisms on the didactic value of the syllabi, on the connections between the various topics (for instance, the fact that the concept of straight line is not connected to that of minimum distance) and on some errors in Euclid.

Amongst the arguments against Euclid there is, in Wilson's opinion, the fact that only England uses Euclid's text; indeed, "France, Germany, Italy and America all agree about geometry teaching as about the teaching of other subjects". The translator of the journal (R. Rubini) notes acidly that: "Thanks to the illuminated Secondary Education Council, Italy has lost this honour, as the council recommends, and thus commands, that Euclid be taught in school !!!" (Wilson 1868). Brioschi and Cremona respond harshly and point out:

*We believe that the logical excellence of Euclid lies in the ordering that you wish to criticise, [...] our schools are intended to produce cultured men and women. We do not aim to teach technical drawing, nor is it important if the students learn this or that theorem, nor that they study many things in a short time. The important thing is that they learn to think, to prove, to deduce. Thus, rapid methods are not useful, and neither are books in which geometry is mixed with arithmetic or algebra. Euclid really is the text which best serves our ends. (Brioschi & Cremona 1869).*

A few years later, Euclid's Elements were no longer considered the required textbook, and, indeed, two calls were made for an elementary textbook, which "keeps rigidly to the Euclidean approach" (Vita 1986). There were ever-growing numbers of "technical" objections to the direct use of Euclid's own text. As far as more general objections are concerned, the analysis of foundations and the presentation of geometry used by various authors (which will see in Hilbert its final realization) no longer permitted, at the turn of the century, the acceptance of definitions and postulates as Euclid presented them. Here, a canonical example would be the tautological Euclidean definition of angle as "mutual inclination of two lines on a plane, which meet one another and are not co-linear". But the problem remained – and remains - as to whether or not, or to what extent, it is right to stick to the Euclidean approach. There are two main considerations:
- How one should define the equality of geometrical figures: in particular, the role that movements or transformations should play in the analysis of geometrical facts.
- Whether or not geometry should be *pure*, i.e., should it be mixed with arithmetic and algebra relying (even if not explicitly) on some prior knowledge of real numbers, or should it be independent.

*4.1.2. Rigid body motions in Italian textbooks at the end of the 19th Century*
In Euclid the theory of the equality of polygons is based, as is known, on the criteria for the equality of triangles, which are in turn based on the concept of *superposition* which makes one figure coincide with another. As is also known, the definition of equality by means of rigid body motion is often circular or makes use of something from outside mathematics. As already mentioned, at the end of the 19th century many well-known Italian geometers started writing geometry textbooks. The ministerial request for a *Euclidean approach* was obeyed. The new Italian texts didn't worry too much about independence of axioms and tried at the same time to maintain some kind of reference to physical space especially in the introduction. Even though they had these common aims, two opposite trends can be seen in these texts. Some authors try to better Euclid as regards rigid body motions, creating new axioms for them. These axioms determine the criteria for the equality of figures. Other authors, more like David Hilbert, prefer to avoid the use of rigid body motions to introduce equality and assume as primitive the concept of equality between certain objects (generally segments and angles), using the concept of 1-1 correspondence.

Of the two approaches, the first seems to have had greater success in schools and is supported by such authors as Achille Sannia and Enrico D'Ovidio (1869), Riccardo De Paolis (De Paolis 1884), Aureliano Faifofer (Faifofer 1878). Supporters of the second approach include Giuseppe Veronese (Veronese 1987) and, to some extent, Enriques and Ugo Amaldi. In fact, one of the most popular school geometry texts written at the turn of the century, Enriques & Amaldi, using rigid body motions to justify the first criterion for the equality of triangles, comments explicitly that they are giving "a simple experimental verification of this criterion and not a proof." (Enriques & Amaldi 1903). The book of De Paolis represents an example of the first approach. His first postulate is about movement (De Paolis 1884):

*Postulate I*
*1. Geometrical figures can move throughout all space*
*2. A figure can move keeping one point fixed*
*3. A figure can move whilst keeping all of its points on a certain line fixed*
*4. To fix a figure it is necessary and sufficient to fix three of its points, not all on the same line.*

De Paolis is able to give, after his first Postulate, the definition of equality, namely:

*1. We will say that two figures are coincident when every point of one is a point of the other and vice versa.*
*2. Two figures are called equal when they can be made to coincide.*

The first criterion (SAS) is then a theorem. De Paolis uses movements, making one pair of sides coincide (the sides AC and A'C' of two given triangles ABC and A'B'C'), and then proceeds in the classical way (that of Euclid), based on rigid body motions.

*4.1.3. The fusion of plane and solid geometry*
As we have seen in postulate I of De Paolis in the previous paragraph, a key methodological question of that period was represented by the so-called "fusionism". The idea originated from projective and descriptive geometry, which worked with projections in space and sections. Foreign textbooks (Bretschneider in 1844, Méray in 1874) had adopted this idea, mixing plane and solid considerations. In Italy, followers of fusionism were, besides de Paolis, Lazzeri & Bassani in 1891. This latter had a great success and was also translated into German by Peter Treutlein in 1911: In this text, the chapter about the properties of incidence refers also to the mutual position of a plane and a straight line, while homothety is defined in space and then in the plane. Pupils were supposed to have a better intuition of spatial relations when passing from space to plane, and to reason by analogy. Moreover, proofs could be presented of plane theorems using configurations in space. For instance, given two triangles with equal bases and heights, the theorem that parallels at the same distance from the base intercept equal segments can be proved using suitable tetrahedra.

## 4.2. Rational arithmetic and real numbers: 1867 – 1923

*4.2.1. Rational Arithmetic*
Rational arithmetic was introduced in high school in 1867, when the Elements of Euclid became the geometry textbook for this school. Rational arithmetic refers to that part of algebra that deals mainly with the properties of integer and rational numbers, exposed in the form of theorems derived from axioms and definitions, in parallel to the rational teaching of geometry (cfr. Menghini, *to appear*). The aim was to present, in the Upper Gymnasium, all mathematics as a deductive science.

In the Technical Institutes the subject was "Arithmetic and Algebra" and the term "rational" didn't appear. But in that period (in 1876 ad 1891) we start to find the term "theory" (of prime numbers, of fractions) in the syllabi. Very likely, there were differences between the Gymnasium and the Technical Institute in the real teaching of this topic, but, in fact, most books of *Rational arithmetic* were devoted both to the first two years of the Technical Institute and to the Gymnasium-Lycée till into the $20^{th}$ century. In fact there was an increase of rigor till to the disappearance of rational arithmetic in 1923 (which would be then limited to the Normal Schools).

"Ambiguous" definitions of number and of their properties could be found in textbooks, as

*- We consider a whole of elements. Each element is a concrete unity. An integer number is a collection of unities*
*- To sum means to add [in the sense of to join] the second number (or quantity) to the first. (Ricotti, 1887).*

Some authors speak of a *collection* or a *class* of objects and of *correspondence* between the objects. So the abstract concept of number arises: In Gigli (1914) the fundamental properties of numbers *a*, *b* and *c* are expressed in words and correspond to:

*- the sum a + b exists*
*- a is equal, minor, or greater than b; a = a; if a = b then b = a;*
*- if a > b then b < a; if a > b, then c exists such that a = b+c*
*- a + b = b + a*
*- (a + b) + c = (a + c) + b*
*[...]*
*- the unity exists; if n is a number n+1 is its successive number.*

The associative property was then a theorem. To prove that $(a+b) + c = a + (b+c)$ following steps were needed

*- (a + b) + c = (b + a) + c (commutative law)*
*- (b + a) + c = (b + c) + a (axiom stated before)*
*- (b + c) + a = a + (b + c) (commutative law)*

In many texts we find also proofs by induction (as in De Franchis, 1912), and in as text of 1908 Catania proposed Peano's arithmetic. It is clear that such a rational approach to arithmetic was at a too high level. Already in 1911, Umberto Scarpis (1911, p.25) wrote:

*Given the strict theoretical development of the teaching of rational arithmetic, the pupils enter the Lycée with very few training in the fundamental rules of calculations, particularly of fractions...*

Many mathematicians asked to abolish this topic and to limit the treatise only to the properties of integers and prime numbers. After the abolishment of rational arithmetic in 1923, Francesco Severi (1931) will explain that algebra could in fact be understood also without rational arithmetic, basing only on the middle school experience.

*4.2.2. The approach to real numbers*
Real numbers were never mentioned explicitly in the programmes of the first two years of the Technical Institutes, nor in those for the physics-mathematics section. We find anyway always mention to square roots and logarithms, sometimes to incommensurables (in 1876 and in the physics-mathematics section in 1885), to limits (1891 and ph.m.s. in 1885), to irrational numbers. In the Lycée real numbers were explicitly mentioned only in 1911. But we find also (in the Lycée and in the first year of the Technical Institute) the theory of *Proportions of magnitudes* (T.I. 1891, Gymnasium 1867 and 1911), needed to treat similarities. The introduction of real numbers by means of the Euclidean approach to proportions and to incommensurables was a typical content for the Lycée, but we could find it also in books for Technical Institutes (as Lazzari & Bassani in 1891).

The programmes for the Technical Institute seemed to suggest a tacit introduction, with mention to square roots and to the way to calculate them (there is for instance a simple introduction in Faifofers's book of Arithmetic, 1883). But in many books of algebra for the Technical Institutes and for the Lycée (Faifofer, 1897; Gigli, 1914/21), we find the complicated approach by means of the so-called *Classi Contigue* (which are not too different from the *nested intervals*). The presentation is mainly geometric and refers to sequences of segments.

Protesting against the curricula proposed in the reform of 1904 (see 2.3) that allowed choosing between Greek and mathematics, mathematics teachers of Mathesis declared that it was scientifically and didactically wrong

*- to treat the equations of the second degree and square roots without having treated irrational numbers*
*- to derive the theory of logarithms from geometric progressions*
*- to treat the theory of measure before that of irrational numbers (Consiglio Direttivo Mathesis, 1909)*

This means that nothing can be done if pupils haven't a complete cognition of the basic topics, in particular of irrational numbers. This point of view was very diffused in Italian school and was the obstacle to an introduction of analytic geometry and of functions at the beginning of high school.

A different position can be found in the suggestions by Castelnuovo for the *Liceo Moderno*, in 1913 (see *2.3.1*). The first three years of this new Lycée, being common to the classical Lycée, presented a rational teaching of arithmetic and geometry. In the last two years a relevant topic was that of approximation, used to introduce incommensurable magnitudes and irrational numbers, but also – in the last year – to introduce definite integrals. As to real numbers, Castelnuovo suggested the following path in the explanations annexed to the programmes (Castelnuovo, 1913): Approximate measures (with measuring instruments); operations with decimal numbers (Castelnuovo invites the teacher to reason about a limited number of numerical examples); comparison between approximate and exact measures; the problem of the existence of a common measure; irrational numbers "about which the teacher will say what is strictly necessary to understand the concept, with only few hints on the operations among them".

But this kind of proposals will not enter in the Italian habits. Still in 1931 Severi (1931) will write that the inquiries of ICMI had shown that in no country real numbers were treated in such a complete way as in Italy. Everywhere only a few hints necessary to applications were given. According to Severi, the introduction of a modern mathematics has to be considered *a priority with respect to methods*.

We will come back on the question of the teaching of real numbers when treating the mathematical contents of the scientific Lycée (section 5).

**4.3 Projective geometry in the Technical Institutes: 1871 – 1923**

In 1871, Minister Castagnola issued new syllabuses for Technical Institutes (Technical Institutes depended at that time on the Ministry of Agriculture, Industry and Commerce) and for the physics-mathematics section (see 2.1). The preface to the syllabi of 1871 shows a deep mathematical knowledge; the mathematician Brioschi was at that time a member of the Council (Consiglio Superiore per l'Istruzione Tecnica) that edited the syllabus, but surely it was inspired by Cremona. It established that mathematics in the Technical Institute should promote "useful and not remote applications" (Ministero etc., 1871), and at the same time "enhance the faculty of reasoning"; and hence the methods used to present it had to be "rigorously precise". Nevertheless, the model was not Euclid's Elements, but was instead provided by the "new doctrine of projectivity", which "supplies graphic constructions to solve first and second order problems" in a straightforward way.

The main points are summarized here: In the first two years of the Technical Institute, common to all the sections, the study of geometry starts from the first elementary notions (angles, circles, inscribed figures, equality, equivalence and similarity of plane figures) and includes the graphical multiplication of segments, the transformation of figures given a scale, preliminary notions of solid figures and their measures. In the third year, geometry for the *physical mathematical section* includes: the theory of projections of geometric forms (projective ranges and pencils, cross ratio, complete quadrilateral) with its applications to the graphical solution of the problems of first and second degree and to the construction of the curves of the second order; the theory of involution; the duality principle in the plane; elements of stereometry and graphic construction of the barycentre of plane figures. In the fourth year the programme provides for the focal properties of conic sections and the projective properties of conics and spheres; the principles of analytic geometry will be founded on the metric relations (by means of the cross ratio) of projective forms. Cremona's *Elementi di Geometria Projettiva*, published in 1873, was written to fulfil this syllabus (see Menghini, 2006).

In the parallel teaching of *descriptive* geometry, the teacher will start from central projection and the projective properties of figures and will treat collinearities, affinities, similarities, with attention to homology, up to the construction of intersections of surfaces of the second degree. Mathematics had six hours per week in the school, descriptive geometry three hours. The only didactical suggestion of the programmes was that pupils had to do a lot of practical work and that the teacher should question them individually and help them in solving exercises.

According to Morpurgo (1875) the syllabi were highly commended. They undoubtedly covered a great deal of ground. Indeed, the original aim was to prepare students in the physics-mathematics classes for direct entry to the School for Engineers (at the third year of the University) without having to attend the first two years in the Faculty of Science. This was in the end not permitted, and in 1876 a new reform took place, based on proposals coming from teacher councils. The physics-mathematics section "preserves its character of school of a general culture, to which the extensive study of Italian letters, that of modern languages, and a strong teaching of the sciences furnish the strength that the humanist education takes from the Greek and Latin literature" (circ. n°119, 1876). The syllabi were reduced. The teaching of projective and descriptive geometry was unified and appeared only in the fourth year. After two years of plane geometry, and one year of solid geometry and trigonometry, the study of projective geometry was whittled down to the study of the projective ranges and pencils, and of the harmonic properties and projective relationships in a circle. Descriptive geometry was restricted to orthogonal and central projections, which were taught together with equalities, similarities, affinities and perspective collinearities.

**4.4. Intuitive geometry in the lower grades from 1881**

*4.4.1. School programmes*
In 1881, intuitive geometry came to life to be taught in the lower Gymnasium (cfr. Menghini, 2010). Previously, geometry was not part of the school programmes for gymnasial students in this age. An earlier intuitive-experimental approach was considered a good help for students to overcome the difficulties they would face when learning rational geometry (see section 4.1). Geometrical drawing should give a further aid. As we read in the decree (R.d., 16-6-1881, n. 323) intuitive geometry had to

*give to youngsters, with easy methods and, as far as possible, with practical proofs, the first and most important notions of geometry, ... useful not only to access geometry, but also to let the students desire to learn, in a rational way, the subject throughout the Lycèe.*

Three years later a decree (R. d., 23-10-1884, n. 2737) of the new minister Coppino abolished the study of intuitive geometry from the lower Gymnasium and moved down rational geometry to the 4th year. This decision was due to the mathematician Eugenio Beltrami. Beltrami's decision (1885) was a consequence of a lack of clear definitions and of the fear that teachers could not emphasize in the right way the experimental-intuitive nature of geometry being tied to the traditional logical-deductive aspect of rational geometry. All this did not touch the *Scuola Tecnica*, where from 1867 a graphic-intuitive method was suggested even to produce simple deductions. The manuals used in this kind of schools always remained of "intuitive geometry", owing the leading ideas to texts as that of Alexis Clairaut, which had been in use in Italian technical schools also before the unification.

In the 1900s a new program was introduced (R. d. 24-10-1900, n. 361) by minister Gallo: intuitive geometry was restored in lower Gymnasium, but, to prevent past problems, the programme included only elementary notions such as the names of the easiest geometrical shapes, the rules to calculate lengths, areas and volumes and also basic geometrical drawing. The guidelines specify that the new studies "were an introduction to rational geometry" and "a review and an expansion of the notions acquired by the pupils at the elementary school"; a practical approach was required, amplified by the teaching of geometrical drawing (but the rules that explain why the various geometrical constructions work had not to be stated).

*4.4.2. The first textbooks of intuitive geometry*
Textbooks of intuitive geometry appeared right after 1900. They present different conceptions. Some authors, as Veronese (1901), introduce – by means of practical examples - the axioms as anticipation to what pupils will see later on, while others, as Giovanni Frattini (1901), present some proofs with practical methods and try to involve the student "imaging" concrete materials. Like all the authors of intuitive geometry books of this period, the straight line is introduced using the idea of a stretched string, a point thinking of grains of sand, a surface thinking of page of a book.

Veronese, after presenting the incidence axioms for plane and space, adds also the congruence axioms:

*Any segment of a line is not congruent to a part of itself, for example the segment AB in the picture is not congruent to*

*CD. This can be verified visually or with a paper strip or compass.*

To avoid infinity, Veronese states that two lines are parallel when they are symmetric about a point, and explains how to verify manually that two lines are parallel. Frattini presents practical proofs in his text, giving also more weight to the properties of polygons. For example:

*There is exactly one perpendicular line through a given point to a line on a plane. Let us bend a plane, imagine an immense piece of paper, and shape right angles so that one folding follows the line we want to draw the perpendicular to, and the other folding must include the point where the perpendicular passes through. Let us reopen the paper; it will be possible to see the trace of the perpendicular through the point and the line (Frattini 1901).*

Or

*The diagonals of a parallelogram bisect each other. Suppose we cut out the parallelogram from a piece of paper, we would have, then, an empty space which could be filled either placing the parallelogram back in the same position or placing the angle A, marked with an arc, on top of the equivalent angle C, the side AD on the equivalent side CB and the side AB on CD. In this way the diagonals of the shape, though upside down, would be in the previous position, the same for their crossing point. The two segments OC and OA would switch their positions: this means they are the same length (Frattini, 1901).*

Veronese only shows diagrams and then lists elementary definitions for triangles, quadrilaterals, and other polygons and for the circle, without stating any property of these shapes.

In many books of intuitive geometry at the beginning of the 20th century we find a very effective use of geometric transformations (isometries). They were considered a suitable tool for an intuitive introduction to geometry. Motions can in fact be carried out experimentally. All authors use them to transport a segment, or to move a figure; only Veronese explicitly defines the symmetry about a point. Both Veronese and Frattini present ruler and compass constructions at the end of the book, as, for instance, the construction of the perpendicular to a straight line passing through a point outside of it. Veronese doesn't justify that construction, while Frattini justifies it referring to the diagonals of a rhombus.

*4.4.3 Further developments*

Some years later a book appeared by Costanzo e Negro (1905). There aren't experimental argumentations nor proofs, but we often find the sentence "the experience teaches and elementary geometry proofs" or "with the usual experimental proof". In 1907, a book by Pisati was published. In the preface the author slightly dissented from the structure of the programmes:

*it seems proved that, in lower middle school, it would be a big mistake to leave the formal aspect of the subject completely apart. Pupils' intellect, in the early years of their life, has a formal nature […]. Certainly, intuitive teaching of geometry is not easier than formal teaching (Pisati, 1907).*

The book presents theorems and classical proofs. We observe that schoolbooks are slowly shifting towards the rigour of rational geometry. Already in 1905, Minister Bianchi had felt the need to remind the teachers in an official note (Boll. Uff. P.I., 1.06.1905) to "avoid abstract statements and proofs" and to use "simple inductive reasoning" to teach the "truths required by the school programmes".

In 1923, the reform made by Minister Gentile turned the clock back. In the first three years of the Gymnasium, geometry studies "must only aim to keep alive all geometrical notions that the pupils have learnt at the primary school and to fix the terminology properly in their memory". Amongst the books published right after the reform of Gentile, we have to mention the textbook by Severi

(1928). Severi's textbook includes a preface by the Minister of Public Education. In spite of the good comments given in the preface, it is difficult to say that the book follows the school programme guidelines. Over the years, middle school geometry had lost its experimental-intuitive nature, and even its terminological function, becoming more and more rational. Textbooks were almost independent from the school programmes. The book by Severi is surely not an exception, although his book for the high school has always been appreciated for the experimental approach (but in fact it included only an experimental approach to axioms). Severi's text for the middle school includes many theorems (also those regarding angles at the centre and angles at the circumference of a circle), with the most traditional proofs, except for using transformations (rotations and symmetries) as a support to the proofs and for avoiding the word "theorem".

In 1940, with the reform of minister Bottai (legge 1-7-1940, n. 899), the first three-years of the Gymnasium, of the Technical School ad of the Complementary School were unified to form the Middle School. With reference to geometry, although its intuitive nature was confirmed, it was suggested to emphasize the evident properties "by means of several suitable examples and exercises, which, sometime, can also assume a demonstrative connotation".

*4.4.4. Ugo Amaldi and Emma Castelnuovo*
An interesting book by Amaldi (1941) followed the reform of 1940. Amaldi had written many geometry textbooks, mainly for the high school, together with Enriques. Amaldi completely stopped the process of "rationalization" of geometry. His textbook is similar to Frattini's book, but it contains some important changes: measurements and geometrical constructions are not illustrated in separate chapters but they are integrated with the other parts of the book, providing a useful didactic tool. We find many figures and references to real life (an opening door gives the idea of infinite planes all passing through the same straight line, paper bands illustrate congruent segments, …), which had completely disappeared in the meantime. So, given the instructions to draw the axis of symmetry of a segment using a ruler and a compass, Amaldi suggests checking the construction by folding the paper and verifying that the circumferences, used for the construction, overlap. To know the sum of the angles of a triangle, he suggests cutting the corners of a triangle drawn on paper, to place them next to each other and to check that they form a straight angle. Similarly, he suggests cutting and folding techniques to verify the properties of quadrilaterals.

The Allied Commission of 1944 (see 3.1.) brought some change in the general organization of Middle Schools (d.m. 14.8.1944), but not many changes in the curricula. The main suggestion, surely coming from Washburne, was to use "active", practical and experimental, teaching methods. These concise programmes didn't influence very much the methods adopted in schools, but an innovation was, indeed, represented in the book on intuitive geometry by Emma Castelnuovo, the daughter of Guido Castelnuovo, in 1948. In her book, E. Castelnuovo follows in Amaldi's footsteps, using drawings, pictures, and cross-references to reality and integration of constructions and measurements. In addition to this, her book interacts with the students, not only to let them follow a logical deduction or a proof, but also to raise questions in their mind.

*What is the meaning – you would question – of the statement that there is only one line passing through two distinct points A, B? How can the contrary be possible? It is true: it is not possible to imagine two or more distinct lines passing through A and B. It is possible, however, to draw with a compass several circles passing through two points […] (Castelnuovo, 1948).*

The book starts with paper folding, and goes on with ruler and square constructions. As Amaldi does, she re-uses the idea of the stretched string to introduce the properties of segments and straight lines; a method already used by Clairaut, who inspired E. Castelnuovo. Simple tools are made up, such as a folding meter to show how to transform a quadrilateral into a different one with the same sides, and to analyze the limit situations.

**4.5 The role of textbooks**

The return to Euclid's Elements, described in 4.1., caused the beginning of an Italian production of manuals. Before that period, besides the text of Legendre already mentioned, we find a text of Franz von Mocnik (mainly in the countries annexed to Austria before the unification), and of Clairaut (in the technical schools of Lombardia). The texts of A. Amiot and Richard Balzer were mainly used in the first two years of the Technical Institutes. The studies on the foundations of mathematics brought in the Italian textbooks a clear tendency towards rigour, as we have seen also in the case of rational arithmetic and even of intuitive geometry. This also depended on the fact that many mathematicians were involved directly in writing textbooks. Klein (1909), though appreciating the scientific level of Italian manuals, strongly criticized them from a pedagogic point of view. In fact textbooks will currently be used by the teachers to follow the theory, and by the pupils only for the exercises. This characteristic will remain in Italian textbooks also in later periods. Textbooks are quite autonomous with respect to the programmes and even "decide" the programmes, as happened particularly in the second half of the 20th century, due to the absence of reforms. Certain traditions of complicated calculations and discussions of parametric equations, particularly frequent in the Scientific Lycée, will be carried out by textbooks and not by programmes.

**5. The reform of 1923 and the birth of the Scientific Lycée**

**5.1. The reform of Giovanni Gentile**

In 1923 the Minister of Education Gentile launched a law on the general reorganization of secondary education. The principles were still those of Casati: the Lycées prepare for the University, while technical instruction is exclusively vocational (see Marchi & Menghini, 2011 and *to appear*).

The physics-mathematics section was suppressed, because of its "hybrid" character and the "discordant duplicity" of its aims (Gentile, 1923), and so was the *Liceo Moderno*, because each school must "correspond to a well defined culture." In their place the *Liceo Scientifico* was established, "in which the needs of scientific culture have their full satisfaction" (Gentile, 1925). In this school mathematics will develop, as Gentile writes, "the systematic and analytical activity of the spirit" (Gentile, 1902).

The *Liceo Scientifico*, unlike other schools, did not have a preparatory course. The access was possible four years after the Elementary School examination, as was for the Technical Institute, which had instead a four-year lower course. The Scientific Lycée had the same subjects of the Classical Lycée (*Liceo Classico*) with the exception of Greek, but it lasted four years as the Technical Institute and shared with this the mathematics programme for the admission to the first class. The only difference, since the Lycée had a more cultural character, was that the admission examination "will require a deeper and more serious mental ability."

The collaboration with the mathematicians in Mathesis was problematic, because they considered as a major reference the programmes of 1918, whose formulation was not in line with Gentile's ideas. Castelnuovo refused to cooperate and his place was taken by Gaetano Scorza.

Scorza's opinions on the role of mathematics in school are antithetical to those of Vailati and

Castelnuovo. "Who understood a proof [...] has the strong conviction of having achieved an absolute truth [...] and of having not obeyed to any other authority than that of his thought." From this point of view that ignores any didactical requirement, mathematics can leave out of consideration "the concrete interpretations of its theories, even [...] those interpretations from which the theories themselves arose" (Scorza, 1921).

The new role of mathematics renders not only possible, but even "necessary [...] a substantial reductions of the syllabi" (Gentile, 1902). In the *Liceo Scientifico* mathematics has 5 hours in the first class and only 3 in the next three classes, compared to 4 hours assigned to Latin for each class. To ensure the freedom of teaching curricula are not subdivided into years, but they provide the arguments of the final state examination.

The mathematical topics for the *Liceo Scientifico* are divided into a first "algorithmic" and a second "theoretical" part: Solving exercises that require only "the direct application of formulas and theorems" will show the knowledge of the notions contained in the first part. The knowledge of the theories contained in the second part, "which are most appropriate to test the candidate's ability to fully understand a rigorous deductive system," will be shown indicating "their general logical development" and exposing "the proofs of the related theorems" (R.d. 06/05/1923). Thus logarithms, defined through arithmetic and geometric progressions, appear, with the use of logarithmic tables, in the algorithmic part. However, the exponential equations, which represent their natural application, are only in the theoretical part, with the real numbers. Formerly, in the *Liceo Moderno*, after introducing powers with a real exponent, logarithms were derived from the exponential equations. In the physics-mathematics section logarithms were defined through the progressions and applied to the calculation of simple and compound interest.

In the algorithmic part we find trigonometry and spherical triangles, which were applied to navigation in the Technical Institutes. Now no reference to theory or applications justifies their presence. Arithmetic (now limited to natural numbers) is treated more extensively than in previous programmes and includes the function of Euler and Fermat's theorem, but it is optional. Rational arithmetic is suppressed, perhaps correctly. Also optional are real numbers, introduced with the "*operations over them*", as happened in the physics-mathematics section.

The programmes of the *Liceo Classico* (which includes the Upper Gymnasium) are substantially a reduction of those of the *Liceo Scientifico*; infinitesimal calculus is not present. Moreover, in the Upper Gymnasium subjects are present that are required for the admission to *Liceo Scientifico* and Technical Institute, as the systems of linear equations and plane Euclidean geometry. Like other schools, even in the *Liceo Classico* there is a distinction between arguments that require only the application of formulas, in particular relating to algebra, and arguments for which the theory has to be exposed, related mainly to Euclidean geometry and proof. The programmes of the Technical Institutes are similar to those of the *Liceo Classico*, but without trigonometry (this is undoubtedly strange) and similarity in space. Mathematics is taught only in the first two years of the new Technical Institute.

The Normal School is replaced by an *Istituto Magistrale*, which is no longer annexed to a Kindergarten or an Elementary School. Work experience is thus missing. The programme mostly overlaps with that of the *Liceo Classico*, but it includes rational arithmetic and does not include real numbers (which are anyway optional in the *Liceo Classico*). Rational arithmetic is considered necessary because the future teacher is to be aware of the laws that underlie the arithmetic that he will teach to the pupils. This subject has remained unchanged even in the successive programmes of this school.

**5.2. Successive reforms till to 1945**

The interest of Fascism in technical and scientific culture led to a change of the mathematics syllabi taking into account those of 1918, as soon as Gentile leaves the government. In 1925 (Minister Fedele) Luigi Bianchi and Scorza, appointed by Gentile in 1923, are still members of the Higher Council. In this year Cartesian coordinates and graphs of functions are introduced in the *Liceo Scientifico*, as formerly in *Liceo Moderno,* independently from differential calculus. Even the physical and mechanical interpretation of some simple elementary functions is suggested before their study with limits and derivatives. All the concepts of calculus are associated with their geometrical or physical meaning. The exponential equations are related to logarithms. The construction of real numbers is now linked to the geometric proportions.

In 1936 (Minister De Vecchi) programmes are edited, which with a few changes have been in effect until 2010. The structure of these programmes is new: they present a subdivision by years and do not contain the distinction between "algorithmic" and "theoretical" topics. Their content is more "practical" than the previous ones: in particular, real numbers are introduced as decimals and elements of probability are present, as it was in the physics-mathematics section. The introduction of real numbers as decimals is motivated by the wish to anticipate the introduction of analytic geometry (which will fail). There is even a reference to the 1918 programmes of this School in the recommendation to use "superior" tools: "So the concepts of integral calculus will allow to find again, in a simple way, rules for calculating areas and volumes, and the knowledge of derivatives will allow to clarify physical concepts and to solve many problems." In 1937 (Minister Bottai) the arguments of the oral state examination for high school are specified, while leaving unchanged the 1936 programs. Their main purpose seems to consist in the suggestion not to abandon the geometric approach in the construction of real numbers by writing "decimal number as a ratio of two quantities."

After WWII, in the period 1944 – 46, the Commission of the Allied Countries (see 3.1. and 4.4.) formulated new syllabi for all the schools. Differently from those of Elementary School, the secondary school syllabi were not promulgated by a law, but were distributed to schools by means of booklets edited by private publishers. Though being the official programmes (till 2010), they are not uniform. For example, only some of them contain the short but interesting premises, which suggest giving way to intuition, to the psychological origin of the concepts, to physical reality. Besides the premises, the syllabi don't present really new elements. The programme of the last four classes of the *Liceo Scientifico* is that of 1936, with minor modifications. Progressions are reduced to "hints"; the distance between skew straight lines is omitted, as are also spherical geometry and probability. The duration of the *Liceo Scientifico* is now five years: the first year coincides with the first two years of the *Liceo Classico,* including subjects that were formerly required for the admission, such as Euclidean plane geometry and systems of linear equations.

In 1952 a large Commission was set up for a general reform of the school (*Commissione Gonnella*). This commission produced interesting proposals (Vita, 1986) but, as had happened with the *Commissione Reale* in 1908, these never were transformed into a law.

**6. Attempts for a reform: 1950 to 2010**

**6.1. Influences of the modern math movement in Italy**

The modern mathematics movements, started in the middle of the 20$^{th}$ century in Europe and in USA, will have an influence on the Italian teaching in Elementary School (see 6.3.2.) and at secondary level. In a first period attention was focused on high school:

*6.1.1. The Royaumont conference.*
In 1958 the O.E.C.E. (European Organization of Economic Cooperation, today O.E.C.D.) organized in Royaumont (Paris) a conference dedicated to the reform of the teaching of mathematics. Each of the 18 member countries was invited to participate with representatives from university and school (the Italian representatives were Luigi Campedelli and E. Castelnuovo). The conference was monopolized by the French School of thought, i.e. by the Bourbakist point of view: Geometry is inserted in the unified organization of mathematics because *it is the study of a vector space of finite dimension on the field of real numbers, eventually provided with a scalar product.* One year later, the same O.E.C.E. followed up the Royaumont conference, by organizing at Dubrovnik (Yugoslavia) a work session in which a group of experts (Mario Villa represented Italy) was charged with elaborating a modern program for the teaching of mathematics in the secondary school (OECE 1962).

*6.1.2. The proposals for reform in Italy.*
The O.E.C.E. conferences, also followed by a conference in Aarhus (Denmark) on the teaching of geometry and by an ICMI conference in Bologna (BUMI 1962), provoked a discrete reform movement in Europe, even at official levels. With regards to Italy, the Ministry instituted in different orders and types of secondary school "pilot classes" for the experimentation of a modern teaching of mathematics. The teachers of these classes could take advantage of special refresher courses and volumes regarding the essential arguments of abstract algebra, logical structure of geometry, and applied mathematics (Villa, 1965).

Despite the push from the other European reforms, the Italian proposals took a different course. Proposals for new programs for the high school followed, in view of a reform that was to take effect in October 1966, and which never took place (in 1962 the Middle School was reformed). The journals *Archimede* and *Periodico di Matematiche* published numerous papers since 1965 (in particular AA.VV., 1965). Programmes were developed in meetings in Gardone and Camaiore; sessions were organized by the Study Centre of the Ministry and the CIIM (Italian Commission for the Teaching of Mathematics, a sub-commission of UMI, the Italian Mathematical Union). Prominent figures in this period were Ugo Morin, vice president of UMI, Campedelli, President of CIIM, Tullio Viola, president of Mathesis, E. Castelnuovo, middle school teacher, Villa, director of the National Education Centre. Also participated Bruno de Finetti, director of the *Periodico di Matematiche*, organ of Mathesis, and Roberto Giannarelli, director of *Archimede*.

The programmes focused primarily on the first two years of high school, which had to be considered unique for all kinds of schools. The scope of the programmes was "to keep pace with the new scientific and didactic needs, but also to update notations, expressions, definitions, to the actual scientific level" (cf. Morin 1965); perhaps it was exactly the new expressions and definitions that rendered the proposals a little bit "hermetic". An ample part relative set theory, order structures and algebraic structures had to *precede* the teaching of geometry. The concept of set seemed the new unifying element of mathematics.

There wasn't, on the contrary, agreement on the moment of introduction, nor on the role to give to the "group of congruences". The teaching of geometry had always been Euclidean and the rigorous introduction of measures was relegated to the last years, derived from proportions and geometrical magnitudes. In the proposals for a reform in the 60s most Italians didn't want to abandon this tradition, rather they agreed in giving up the traditional axiomatic setting. For this aim it seemed useful to "recover" Klein's Erlanger Program (not highly considered in the Royaumont conference).

In the years 1966 and 1967 a large commission meeting at Villa Falconieri in Frascati under the auspices of UMI and CIIM ended the long series of meetings dedicated to the programmes with a final proposal (R.G. 1966 and 1967). In the *Frascati programmes* for the last three years of high school there are, among optional topics, the "elementary transformations and their groups", and even the "projective enlargement of the affine or Euclidean space", which enter in the study of the geometric vector plane in the 3rd year, and in its extension to space with scalar product in the 4th year. Many reports of the teachers of the pilot classes, or at least of the experts of didactics (cfr. Castelnuovo, 1964; Mancini Proia, 1967) sustain that the road followed is the one suggested by Klein, who is placed among the innovators of the 60's; and there is also explicit reference to the Belgian Paul Libois, to whom they owed several models for presenting geometric transformations (Libois, 1963). There are also references to Georges Papy and Gustave Choquet. At the end of the 60's the concerns for a new axiomatic system made one appreciate the work of Choquet, who presented, in the appendix of his book on the teaching of geometry (1964), a system of axioms for the metric plane, based on distance and symmetry axioms.

It must be emphasized that neither the many articles published in journals, nor the work of the 82 pilot classes in mathematics, were sufficient to spread and make the new ideas appreciated among the majority of mathematics teachers. In the absence of a reform, the more "Bourbakist" ideas slowly took a back position. Klein and Choquet survived and we find traces of both in the successive proposals and in the conception of many new books.

## 6.2 Experimental programmes

In some Italian cities, University – School groups were formed for the testing of new curricula in the high school. These groups became official starting from the years '75-'76 with a contract between UMI and CNR (*Consiglio Nazionale delle Ricerche*). There was also a contract between CNR and Mathesis (OPI and RICME projects for the Elementary School), and a contract with the University of Genoa for the Middle School.

More than modern mathematics, methodology seemed important as well as the freedom to experiment new topics; teaching via problems appeared particularly effective. These projects proceeded independently, albeit with many opportunities for discussion and with common topics. They tested their proposals asking, if necessary, a permission from the Ministry. They produced materials and textbooks, but did not lead to a unified proposal for a reform, and reached only a small part of the Italian teachers.

For a comprehensive and widespread proposal it is necessary to wait until 1985. In that year Minister Falcucci promotes the PNI (Piano Nazionale per l'Informatica), which – with the official aim to promote the teaching of computer science in school – proposed new syllabi for mathematics and physics to allow the use of a new technology (Ministero Pubblica Istruzione, 1991 e 1992). The programmes of PNI were experimental, but they were strongly sustained by the Minister and in most Italian schools there were classes with the new programmes. Moreover many teachers were trained to teach in these classes; this was the first *national in-service teacher training* in Italy. The PNI programmes lasted till to the reform of 2010 (see also Linati, 2011).
They represent the real novelty after the Gentile reform. Among the contents we find topics of the programmes of Frascati, and of the many Italian projects. The teaching of functions and of the Cartesian plane is precocious, and furnishes a geometric interpretation for the solutions of equations. Euclidean geometry can be taught by means of geometrical transformations. Elements of probability and statistics enter the curriculum, and of course informatics, whose algorithmic aspects were considered particularly relevant. But slowly technologies as an instrument to teach prevailed. Important general ideas were those of the mathematization of real phenomena and of the teaching via problems (Ciarrapico, 2002).

In the same period the absence of school reforms brought the Minister to form a Commission, guided by Beniamino Brocca, in order to formulate a general reform of secondary instruction (Middle School had been in the meantime re-reformed in 1979). As to mathematics, the proposed programmes were very similar to those of PNI (many mathematicians and expert teachers were in both commissions); we underline only the addition, in the last three years, of elements of arithmetic, which was considered a rich field of problems.

## 6.3 Reforms of compulsory and technical schools

As we have seen in section 6.2., there have been no laws, after 1936, concerning the programmes of high school. In 2010 a general reform of High School (*Indicazioni Nazionali*) has been promulgated (archivio.pubblica.istruzione.it/riforma_superiori/nuovesuperiori/index.html). Other school types had been instead reformed in previous years:

### 6.3.1. Middle school
An important reform, after the Bottai reform of 1940 (section 4.4.3) and the reform extending Middle School to *all* pupils in 1962 (*Scuola Media Unica*), was the one of 1979 (d.m. 9.2.79). As in 1962 (see 3.4), mathematics is taught together with the experimental sciences, for a total of 6 hours per week. The programmes were new, and they showed that the Bourbakist wave was definitely ceased (Ciarrapico, 2002). Among the aims for the teaching of mathematics there is the development of intuition, and the skill to communicate using verbal, graphical and symbolic expressions. Teaching must start from concrete facts and proceed via problems. Contents are divided into seven themes: Geometry and the physical world, sets of numbers, mathematics of certain and probable events, problems and equations, the method of coordinates, geometrical transformations, correspondences and structural analogies.

Computer Science isn't among the contents, but there is a reference to the use of calculating instruments. The programmes were surely inspired by the work of E. Castelnuovo, who was in the commission. They were very vast and this – together with the new topics - will represent a problem. Many teachers weren't able to fulfil them and mostly continued to teach traditional algebra and geometry.

*6.3.2. Elementary School*
Till to a reform of 1955, the mathematical contents of Elementary School had been chosen thinking at their use in life (simple calculations, measures), considering that in that period many pupils ended the school after Elementary School, even if three more years were compulsory. In the successive reform of 1985 (D.P.R. 105/85) the conception of mathematics teaching had completely changed. These programmes were also inspired by research in mathematics education, so the need was recognized for a long path, which, starting from the observation of reality, brings to mathematization, resolution of problems, till to mathematical abstraction. We find five themes: Problems – arithmetic – geometry and measure – logic – probability, statistics and informatics. The themes "Problems" and "Logic" don't have proper contents, but show a methodology. Before these programmes, and in the absence of a proper reform, many teachers and textbooks had adopted a Bourbakist point of view introducing (much more than in secondary school) the concepts of set theory in order to present the cardinal concept of number. The programmes of 1985 tried to limit this habit by writing that "the formal symbolism for the logical –set-theoretical operations is not necessary to introduce the integer number and the operations among them".

*6.3.3. Technical Institutes*
It is difficult to give an account of the programmes of the Technical Institutes. The number of their sections raised after 1923, the teaching of mathematics covered from two to five years and among traditional arguments we can find, in different sections, complex numbers (industrial and nautical section); descriptive geometry (industrial section); financial and actuarial mathematics, combinatorics and probability (commercial section); surfaces of the second degree, helices and helicoids (section for geometers); spherical trigonometry and the sphere (nautical section). In all the sections the teaching of geometry shouldn't be strictly rational.

The Italian law establishes that the Minister, without a law of the Parliament, can change the programmes for the Technical Institutes. This brought to a larger number of changes with respect to the Lycée. We don't assist to great revolutions, but – being the programmes of different sections prepared by different commissions – it can happen, as in 1966 for the new section of "correspondent in foreign languages", that we find really interesting programmes. So the D.M. 8-8-1966 insists on the necessity to link the teaching of mathematics to that of other disciplines. "It will be useful", we read, "starting from a concrete experience, to give the habit to reason [...] to generalize". Algebra is then to be treated in parallel with the analytic representation in the Cartesian plane. Also, the historical development of the theories has to be considered.

Something new happened in the 70ties, when the computer widened its influence and a new section oriented to computer science appeared. In this section computer science is a subject for its own, but the parallel teaching of mathematics is strengthened and contains statistics, probability and linear programming. Other interesting syllabi for the Technical Institutes will appear after the 80s, with many topics in common with the PNI. Moreover Technical Institutes adopted the norms for the *minimal mathematical contents* in compulsory school, i.e. till to 16 years, established by the Ministry in 2007 (*Assi culturali*, d. 22.08.2007, n.139).

# Acknowledgements.
Many thanks to Claudio Bernardi and Fulvia Furinghetti for their collaboration and their careful review.

# Glossary

**Arithmetic:**

**Combinatorics:** It concentrates on counting the number of certain combinatorial objects. It provides a framework for counting permutations, combinations and partitions.

**Fermat's theorem:** In number theory, Fermat's Theorem states that no three positive integers a, b, and c can satisfy the equation $a^n + b^n = c^n$ for any integer value of n greater than two.

**Function of Euler:** In number theory, Euler's phi function $\varphi(n)$ is an arithmetic function that counts the number of positive integers less than or equal to n that are relatively prime to n.

**Peano's arithmetic:** In mathematical logic, the Peano axioms are a set of axioms for the natural numbers. These axioms are used for a formal introduction of arithmetic.

**Bourbaki:** Nicolas Bourbaki is the collective pseudonym under which a group of (mainly French) 20th-century mathematicians wrote a series of books presenting an exposition of modern advanced mathematics, beginning in 1935. With the goal of founding all of mathematics on set theory, the group strove for rigour and generality.

**Conic sections and curves of the second order:** In mathematics, a conic section (a conic) is a curve obtained as the intersection of a cone with a plane. In analytic geometry, a conic may be defined as a plane algebraic curve of degree 2. A conic consists of those points whose distances to some point, called a focus, and some line, called a directrix, are in a fixed ratio, called the eccentricity. The three types of conic section are the hyperbola, the parabola, and the ellipse.

**John Dewey:** (October 20, 1859 – June 1, 1952) was an American philosopher, psychologist and educational reformer whose ideas have been influential in education and social reform. Dewey was one of the founders of functional psychology. He was a major representative of progressive education and liberalism.

**Erlanger Program:** Presented by Felix Klein at Erlangen in 1872, it proposes a unified view of geometry as the study of the properties of a space that are invariant under a given group of transformations,

**Functions and rate of growth**: In mathematics, a function is a relation between a set of inputs and a set of potential outputs with the property that each input is related to exactly one output. The output of the function f corresponding to an input x is denoted by f(x). The rate of growth of a function in a given interval [a,b] corresponds to its slope: (f(b) – f(a))/(b-a). The rate of growth of a function in a given point is a limit process when b tends to a and is calculated by means of the derivatives.

**Helices and helicoids:** A helix is a type of smooth space curve, i.e. a curve in three-dimensional space. A "filled-in" helix – for example, a spiral ramp – is called a helicoid.

**Infinitesimal calculus:**

**Limits:** In mathematics, the concept of a "limit" is used to describe the value that a function or sequence "approaches" as the input or index approaches some value. Limits are essential to calculus and are used to define continuity, derivatives, and integrals.

**Derivatives and integral calculus:** Differential calculus (derivatives) is concerned with the study of the rates at which quantities change; integral calculus is concerned with the calculatuon of the area of the region in the xy-plane bounded by the graph of a function f, the x-axis, and the vertical lines x = a and x = b (definite integral).

**Maria Montessori**: (August 31, 1870 – May 6, 1952) was an Italian physician and educator, best known for the philosophy of education which bears her name. Her educational method is in use today in public and private schools throughout the world.

**Projective geometry:** In mathematics, projective geometry is the study of geometric properties that are invariant under projective transformations. One source for projective geometry was indeed the theory of perspective. Its study includes:

**Complete quadrilateral:** Dually, a complete quadrilateral is a system of four lines, no three of which pass through the same point, and the six points of intersection of these lines. It is the projective *dua*l of a complete quadrangle is a complete quadrilateral. The four points on the line deriving from the sides and diagonals of the quadrangle are called a *harmonic range*.

**Cross ratio.** In geometry, the cross-ratio, also called double ratio and anharmonic ratio, is a special number associated with an ordered quadruple of collinear points, particularly points on a projective line.

**Graphical multiplication of segments:** Multiplication and division are performed by using proportional segments, which are cut off by parallel lines on the sides of an angle.

**Harmonic division**. In geometry, harmonic division of a line segment AB means identifying two points C and D such that AB is divided internally and externally in the same ratio. It is also related to the cross-ratio: Given this metric context, the harmonic range is characterized by a cross-ratio of minus one:

**Involution:** An involution is a projectivity of period 2, that is, a projectivity that interchanges pairs of points.

**Polarity:** In geometry, the terms pole and polar are used to describe a point and a line that have a unique reciprocal relationship with respect to a given conic section. If the point lies on the conic section, its polar is the tangent line to the conic section at that point.

**Realschule:** type of modern or technical secondary school in Germany, Austria, and other countries, which offers a more scientific curriculum compared to the Gymnasium.

**Spherical trigonometry:** A branch of spherical geometry which deals with polygons (especially *triangles*) on the sphere and the relationships between the sides and the angles. It is of great importance for calculations in astronomy and earth-surface, orbital and space navigation.

## Biographical Sketch

**Marta Menghini** received her degree in Mathematics from the University of Rome "La Sapienza" on July 1976. She is associate professor in the Department of Mathematics of the same University. She produced papers in the field of Combinatoric Geometry and of History of Mathematics (in particular on the development of projective and algebraic geometry in the last century). She was a member of the commission of the Italian National Seminar of Didactics of Mathematics; from 1987 till 1995 she was in the Editorial board of "Epsilon", a magazine for high school teachers of scientific disciplines.
From 1987 till 2000 she was Scientific Director of a research project about Innovation of Mathematics at School (due to an agreement between CNR and University) and edited many works of the teachers involved in the project (teaching materials, reports, conferences). From 1998 to 2000 she took part to the European Project "Modem", about innovation and production of materials for didactics of mathematics. She was in the Commission of the Italian Mathematical Union for the curricular reform of the High School.
In 2008 she was in the scientific and in the organizing committee of the Congress held in Rome in the occasion of the centannial of ICMI.
She is the author of numerous published works on mathematics education in secondary school and on the history of mathematics education, in particular on the teaching of geometry and on the history of teaching geometry.